\documentclass [11pt,a4paper, reqno]{article}
\usepackage{subfigure, graphicx, tabularx}
\usepackage{amsfonts}
\usepackage{amsmath}
\usepackage{amssymb}
\usepackage[english]{babel}
\usepackage{amsfonts}
\usepackage{array, natbib}
\makeatletter
\def\hlinewd#1{%
\noalign{\ifnum0=`}\fi\hrule \@height #1 %
\futurelet\reserved@a\@xhline}
\makeatother
\newlength{\abstractwidth}
\newtheorem{theo}{Theorem}[section]

\newtheorem{prop}{Proposition}[section]
\newtheorem{cor}{Corollary}[section]

\renewcommand{\baselinestretch}{1.31}
\begin{document}

\def\Cov{\mbox{Cov}}
\def\N{{\rm I\!N}}
\def\E{{\rm I\!E}}
\def\ds{\displaystyle}
\def\ind{ {{\rm 1}\hskip-2.2pt{\rm l}}}
\def\att{ { {\huge $\ds{\triangle}$} \hskip-14pt {\large !} } }
\def\a{\alpha}
\def\b{\beta}
\def\e{\varepsilon}
\def\d{\delta}
\def\D{\Delta}
\def\f{\phi}
\def\h{\eta}
\def\L{\Lambda}
\def\l{\lambda}
\def\p{\rho}
\def\oo{\infty}
\def\||{\parallel}
\def\R{\mathbb{R}}
\def\P{\mathbb{P}}
\def\tP{\widetilde{P}}
\def\u{\mathbf{u}}
\def\x{\mathbf{x}}
\def\z{\mathbf{z}}
\def\Z{\textbf{Z}}
\def\X{\mathbf{X}}
\def\W{\textbf{W}}
\def\w{\textbf{w}}
\def\vv{\mathbf{v}}
\def\({\left(}
\def\){\right)}
\def\1{\mathbf{1}}
\def\Qn{\bar{Q}_n}
\def\Qnn{\widetilde{Q}_n}
\def\Q{\widehat{Q}_n}
\def\Unk{U\(\frac{n}{k}\)}
\def\Uunk{\widehat{U}\(\frac{n}{k}\)}
\def\aa{\widehat{\a}}
\def\ppsi{\widehat{\psi}_n}
\def\psis{\psi_n^*}
\def\psij{\psi_{n,j}^*}
\newcommand{\gr}{ \gamma \in I\!R }
\newcommand{\bb}{ {\bf b} }
\newcommand{\bee}{ {\bf e} }
\newcommand{\nn}{ n \rightarrow \infty }
\newcommand{\no}{ \noindent }
\newcommand{\xx}{ x \rightarrow \infty }
\newcommand{\ks}{ \vspace{0.3cm} \noindent }
\newcommand{\bY}{ {\bf Y } }
\newcommand{\bm}{ \mbox{\boldmath $\mu$} }
\newcommand{\bt}{ \mbox{\boldmath $\tau$} }
\newcommand{\be}{ \mbox{\boldmath $\epsilon$} }
\newcommand{\tbe}{ \mbox{\boldmath $\tilde{\epsilon}$} }
\newcommand{\bal}{ \mbox{\boldmath $\alpha$} }
\newcommand{\bRk}{ {\bf R} }
\newcommand{\bC}{ {\bf C} }
\newcommand{\bc}{ {\bf c} }
\newcommand{\bD}{ {\bf D} }
\newcommand{\bon}{ {\bf 1} }
\newcommand{\bff}{ {\bf f} }
\newcommand{\tbg}{ \tilde{{\bf g }} }
\newcommand{\bZ}{ {\bf Z} }
\newcommand{\bX}{ {\bf X} }
\newcommand{\argmax}{ \mbox{argmax} }

\newcommand{\hr}{\hat{\rho} }
\newcommand{\bS}{ {\bf S}^{-1/2}_{k} }
\newcommand{\ba}{ {\bf a}_{kn} }
\newcommand{\tbX}{ \tilde{\bX} }
\newcommand{\htr}{ \hat{\tilde{\rho}} }
\newcommand{\tr}{ \tilde{\rho} }
\newcommand{\tee}{ \tilde{\epsilon} }
\newcommand{\ee}{ \varepsilon }
\newcommand{\ccrit}{c_{crit}}

\title{\bf {\LARGE  Improved estimation of the extreme value index \\ using  related variables}}
\author{\Large {Hanan Ahmed}\\{\large{\it Tilburg University}} \\ \and {\Large John H.J.\ Einmahl}\\
{\large{\it  Tilburg University}}\\
\\
}

\maketitle

{\normalsize{\bf Abstract}. 
Heavy tailed phenomena are naturally analyzed by extreme value statistics. A crucial step in such an analysis is the estimation of the extreme value index, which describes the tail heaviness of the underlying probability distribution. We consider the situation where we have next to the $n$ observations of interest another $n+m$ observations of one or more related variables, like, e.g., financial losses due to earthquakes and   the related amounts of energy released,  for a longer period than that of the losses. Based on such  a data set, we present an adapted version of the Hill estimator that shows greatly improved behavior and we
establish the  asymptotic normality of this estimator. For this adaptation the tail dependence between the variable of interest and the related variable(s) plays an important role. A    simulation study confirms the substantially improved performance   of our adapted estimator relative to the Hill estimator. We also present an application to  the aforementioned earthquake losses.}
\bigskip

\noindent Keywords: Asymptotic normality, Heavy tail,  Hill estimator,     Tail dependence,     Variance reduction 

\newpage   
           
\section{Introduction}
%%%%%%%%definitions%%%%%%%%%%
Consider univariate extreme value theory for heavy tails, that is,  the case where the extreme value index $\gamma$ is positive. This index describes the tail heaviness of the underlying probability distribution, the larger $\gamma$, the heavier the tail. See \cite{de2006extreme} or \cite{beir} for a comprehensive  introduction to univariate and multivariate extreme value theory. Given a random sample, we estimate $\gamma$ with the well-known and often used \cite{hill} estimator. Such an estimate of $\gamma$ is the  crucial ingredient for estimating important tail functionals of the distribution, like very high quantiles, very small tail probabilities, but also the Expected Shortfall or an excess-of-loss reinsurance premium.

In this paper we consider the situation where we have a bivariate (or multivariate) data set, where the first component is the variable of interest with extreme value index $\gamma_1$ and the second component is a heavy-tailed  related variable, with extreme value index $\gamma_2$, that should help to improve the estimation  of $\gamma_1$. We assume that we have a random sample of size $n$  of these pairs and in addition another $m$, independent of the pairs and mutually independent, observations of the second component. Hence we have a larger sample of the related variable than that of the variable of interest. Such a situation can occur when we have recorded both variables for a certain period of time (2008-2017, say), but in addition have  data for  the second variable only, for an earlier period (1980-2007, say). Here we can think of financial losses as the variable of interest and some physical quantity (like water height, wind speed, earthquake magnitude) as the related variable. Another situation where our setup  can occur is when in a certain period the related variable is measured more frequently than the variable of interest, but other situations might also lead to our setting. We will also consider the setting where there is \textit{more} than one related variable, the multivariate case, but in this introductory section we will focus on the bivariate case.

We can estimate $\gamma_1$ with the Hill estimator $\hat \gamma_1$ and $\gamma_2$ with the Hill estimators $\hat \gamma_2$, based on the $n$ data,  and  $\hat \gamma_{2+}$, based on all $n+m$ data. The latter estimator is  better than $\hat \gamma_2$, ``hence" their difference can be used to update and improve $\hat \gamma_1$. For this updating the strength of the  tail dependence between both variables is important and should be estimated. A detailed derivation of our adapted Hill estimator is presented in the next section. We will show that our estimator improves greatly on the Hill estimator. To the best of our knowledge this approach is novel and there are no results of this type in the literature.

The remainder of this paper is organized as follows. In Section 2, for the clearness of the exposition,  the  bivariate case is treated as indicated above and the asymptotic normality of the adapted estimator is established and in Section 3 the corresponding results for the  multivariate case are presented.  In Section 4, the finite sample performance of our estimator is studied through a simulation study, which  confirms the improved performance of the adapted Hill estimator. In Section 5, we present an application to   earthquake damage amounts with  the ``amount of energy released" as related variable. The proofs of the results in Section 3 are deferred to Section 6.  Since Section 3 generalizes Section 2, the proofs of Section 2 can be obtained by specializing  those of Section 3, and are hence omitted.

\section{Main results: the bivariate case}
Let $F$ be a bivariate distribution function with marginals $F_1$ and $F_2$. Assume that $F$ is in the bivariate max-domain of attraction (i.e., $F \in D(G)$) with both extreme value indices $\gamma_1$ and $\gamma_2$ \textit{positive}, see Chapter 6 in \cite{de2006extreme}. Let $U_j=F^{-1}_j(1-1/\cdot)$ be the tail quantile corresponding to $F_j, j=1,2$. Then $F \in D(G)$ with positive extreme value indices implies that $U_j$ is regularly varying with index $\gamma_j$, $j=1,2$, i.e., $\lim_{t\to\infty} U_j(tx)/U_j(t)=x^{\gamma_j}, x>0.$   Let $(X_1,Y_1)$ have distribution function $F$. Then $F \in D(G)$ also implies the existence of the tail copula $R$ defined by
\begin{equation}\label{tcop}
R(x,y)=\lim\limits_{t\downarrow 0}\frac{1}{t} \mathbb{P} \left(1-F_1\left(X_1 \right) \leq tx,1-F_2\left(Y_1 \right) \leq ty \right) \quad  (x,y)\in [0,\infty]^2\setminus \{(\infty, \infty)\}.\end{equation}

Let    $(X_1,Y_1),\ldots,(X_n,Y_n)$ be   a bivariate random sample from $F$,  and let $Y_{n+1},\ldots,Y_{n+m}$ be a univariate random sample from $F_2$, independent from the $n$ pairs. Denote the order statistics of the $X_i, i=1, \ldots, n$, with $X_{1,n}\leq \ldots \leq X_{n,n}$ and   use similar notation for the order statistics of  the  $Y_i, i=1, \ldots, n$, and also for the order statistics of all the $Y_i, i=1, \ldots, n+m$.  For  $k\in \{1, \ldots, n-1\}$ define the \cite{hill} estimator of $\gamma_1$ by
\begin{equation}\label{hil}\hat\gamma_1=\frac{1}{k}\sum_{i=0}^{k-1}\log X_{n-i,n}-\log X_{n-k,n}.\end{equation}
Define, using the same $k$, similarly  the  Hill estimator $\hat\gamma_2$ based on the $Y_i, i=1, \ldots, n$, and let  $\hat\gamma_{2+}$ be the Hill estimator of all  $Y_i, i=1, \ldots, n+m$, with $k$ replaced by $k_+\in \{k+1, \ldots, n+m\}$.
Throughout for the asymptotical theory we will assume that $m=m(n)$ and that
\begin{equation}\label{kaa}
k\to \infty, \,\frac{k}{n} \to 0,\,  \frac{k_+}{n+m}\to 0, \, \sqrt{\frac{k}{k_+}}\to \nu\in (0,1),\,  \frac{n}{n+m}\frac{k_+}{k} \to\beta\in (0,1], \quad \mbox{as } n \to \infty.
\end{equation}
Observe that we now also have $k_+\to \infty$ and $m\to \infty$; actually $n/(n+m)\to \beta \nu^2\in (0,1)$.

First we consider the joint
  asymptotic normality of the three Hill estimators $\hat\gamma_1, \hat\gamma_2$, and  $\hat\gamma_{2+}$. For this, we need the usual second order conditions, on $F_1$ and $F_2$: there exist positive or negative functions $A_j, j=1, 2,$ with $\lim_{t\to\infty}A_j(t)=0$,  such that for $x>0$
\begin{equation}\label{so}
    \lim_{t\to\infty}\frac{\frac{U_j(tx)}{U_j(t)}-x^{\gamma_j}}{A_j(t)}= x^{\gamma_j} \frac{x^{\rho_j}-1}{\rho_j}, \quad \mbox{for some } \rho_j\leq 0,\,  j=1,2.
 \end{equation}

\begin{prop}
If $F\in D(G)$,  conditions (\ref{kaa}) and (\ref{so}) hold,  and $ \sqrt{k}A_j(\frac{n}{k})\to\lambda_j \in \mathbb{R}, j=1,2$, %and $ \sqrt{k_+}A_2(\frac{n+m}{k_+})\to\lambda_{2+}\in \mathbb{R}$,
as $n\to\infty$, then
\begin{equation}
\left(\sqrt{k}(\hat{\gamma}_1-\gamma_1),\sqrt{k}(\hat{\gamma}_2-\gamma_2),\sqrt{k_+}(\hat{\gamma}_{2+}-\gamma_{2})\right)\xrightarrow{d} N\left ( \left(\frac{\lambda_1}{1-\rho_1}, \frac{\lambda_2}{1-\rho_2} ,\frac{\lambda_2\beta^{-\rho_2}}{\nu(1-\rho_{2})} \right),\breve\Sigma\right ), \label{eq1}
\end{equation}
with  $$ \breve \Sigma= \left[ \begin{array}{ccc} \gamma_1^2 & R(1,1)\gamma_1\gamma_2  & \nu R(1,\beta)\gamma_1 \gamma_2  \\ \\ R(1,1)\gamma_1\gamma_2 &  \gamma_2^2 & \nu \beta\gamma_2^2 \\  \\ \nu R(1,\beta)\gamma_1\gamma_2 &  \nu \beta\gamma_2^2 &  \gamma_2^2\end{array}\right].$$
\label{prop1}
   \end{prop}
   \begin{cor}
Under the conditions of Proposition \ref{prop1}, as $n\to\infty$,
\begin{eqnarray*}  &&\left(\sqrt{k}(\hat{\gamma}_1-\gamma_1),\sqrt{k}(\hat{\gamma}_{2+}-\hat{\gamma}_2)\right)\\
&&\qquad \xrightarrow{d} N\left( \left(\frac{\lambda_1}{1-\rho_1}, \frac{\lambda_2(\beta^{-\rho_2}-1) }{1-\rho_2} \right), \left[ \begin{array}{cc} \gamma_1^2 & \big(\nu^2 R(1,\beta)-R(1,1) \big) \gamma_1 \gamma_2 \\ \big(\nu^2 R(1,\beta)-R(1,1) \big) \gamma_1 \gamma_2 &  \big(1+\nu^2-2\nu^2\beta\big)\gamma_2^2 \end{array}\right]\right).\end{eqnarray*}  \label{cor1}
   \end{cor}

Corollary \ref{cor1} is the basis for finding our adapted Hill estimator. Take $\lambda_1=\lambda_2=0$. %and assume that the limiting covariance matrix is invertible.
The tail copula $R$  is estimated as usual,   cf.\ \cite{drees}, by
\begin{equation}\label{rr} \hat{R}(x,y)=\frac{1}{k}\sum\limits_{i=1}^{n}1_{[X_i\geq X_{n-[kx]+1,n}, Y_i\geq Y_{n-[ky]+1,n}]},\quad  x,y\geq 0.\end{equation}
 Now consider $(\hat{\gamma}_1,\hat{\gamma}_{2+}-\hat{\gamma}_2)$ and  its approximate  bivariate  normal distribution according to Corollary \ref{cor1}, with estimated covariance matrix:
  \begin{equation}\label{limnor} N\left ((\gamma_1,0), \frac{1}{k}\left[ \begin{array}{cc} \hat{\gamma}_1^2 & (\frac{k}{k_+}\hat{R}(1,\frac{k_+}{k}\frac{n}{n+m})-\hat{R}(1,1)) \hat{\gamma}_1 \hat{\gamma}_{2+} \\(\frac{k}{k_+}\hat{R}(1,\frac{k_+}{k}\frac{n}{n+m})-\hat{R}(1,1)) \hat{\gamma}_1 \hat{\gamma}_{2+} &  \big(1+\frac{k}{k_+}-2\frac{n}{n+m}\big)\hat{\gamma}_{2+}^2 \end{array}\right]\right). \end{equation}
%$$f(\tilde{\Gamma},\tilde{\mu},\tilde{\Sigma})=\frac{1}{\sqrt{2\pi^3 |\tilde{\Sigma}|}} \exp\left[\frac{-1}{2} (\tilde{\Gamma}-\tilde{\mu})^T \tilde{\Sigma}^{-1} (\tilde{\Gamma}-\tilde{\mu}) \right],$$  \\  where $$\tilde{\Gamma}= \left[\begin{array}{c} \hat{\gamma}_1 \\ \\ \hat{\gamma}_{2+}-\hat{\gamma}_2 \end{array}\right], \quad \tilde{\mu}= \left[\begin{array}{c} \gamma_1 \\ \\ 0 \end{array}\right], \quad \tilde{\Sigma}=\frac{1}{k}\left[ \begin{array}{cc} \hat{\gamma}_1^2 & \tau \hat{\gamma}_1 \hat{\gamma}_{2+} \\ \tau \hat{\gamma}_1 \hat{\gamma}_{2+} &  \big(1+\frac{k}{k_+}-2\frac{n}{n+m}\big)\hat{\gamma}_{2+}^2 \end{array}\right],$$ $\tau=\frac{k}{k_+}\hat{R}(1,\tilde{\beta})-\hat{R}(1,1)$ and $\tilde{\beta}=\frac{k_+}{k}\frac{n}{n+m}$.
Maximizing the corresponding likelihood with respect to $\gamma_1$, we obtain our adapted estimator  for $\gamma_1$:
%$$ \frac{\partial \log f}{\partial \gamma_1}= \frac{-2k}{\hat{\gamma}_1^2}(\hat{\gamma}_1-\gamma_1)+\frac{2 k \tau }{\hat{\gamma}_1\hat{\gamma}_{2+} [1+\frac{k}{k_+} -2\frac{n}{n+m}]}(\hat{\gamma}_{2+}-\hat{\gamma}_2)=0\\ $$
%$$\frac{\hat{\gamma}_1-\gamma_1}{\hat{\gamma}_1}=\frac{\tau }{\hat{\gamma}_{2+} [1+\frac{k}{k_+} -2\frac{n}{n+m}]}(\hat{\gamma}_{2+}-\hat{\gamma}_2)\\$$
%$$\hat{\gamma}_1-\gamma_1=\frac{\hat{\gamma}_1}{\hat{\gamma}_{2+}}\frac{\tau}{[1+\frac{k}{k_+}-2\frac{n}{n+m}]}(\hat{\gamma}_{2+}-\hat{\gamma}_2)\\$$
%$$\gamma_1=\hat{\gamma}_1-\frac{\hat{\gamma}_1}{\hat{\gamma}_{2+}}\frac{ \tau}{ [1+\frac{k}{k_+} -2\frac{n}{n+m}]}(\hat{\gamma}_{2+}-\hat{\gamma}_2)\\$$
%$$\gamma_1=\hat{\gamma}_1+\frac{\hat{\gamma}_1}{\hat{\gamma}_{2+}}\left(\frac{\hat{R}(1,1)-\frac{k}{k+} \hat{R}(1,\tilde{\beta})}{1+\frac{k}{k_+} -2\frac{n}{n+m}}\right)(\hat{\gamma}_{2+}-\hat{\gamma}{_2}).$$
 \begin{equation}
  \hat{\gamma}_{1,2}=\hat{\gamma}_1+\frac{\hat{\gamma}_1}{\hat{\gamma}_{2+}}\left(\frac{\hat{R}(1,1)-\frac{k}{k_+} \hat{R}(1,\frac{k_+}{k}\frac{n}{n+m})}{1+\frac{k}{k_+} -2\frac{n}{n+m}}\right)(\hat{\gamma}_{2+}-\hat{\gamma}{_2}).
  \label{eq2}
  \end{equation}
  The main result of this section, the asymptotic normality of this estimator, shows that it improves substantially on the Hill estimator.

 \begin{theo}
  Under the conditions of Proposition \ref{prop1}, as $n\to\infty$, %and assuming that the limiting covariance matrix in Corollary \ref{cor1} is invertible, we have
 $$ \sqrt{k}(\hat{\gamma}_{1,2}-\gamma_1)\xrightarrow{d}N\left(\frac{\lambda_1}{1-\rho_1}+\frac{\gamma_1}{\gamma_2}\cdot  \frac{R(1,1)-\nu^2 R(1,\beta)}{1+\nu^2 -2\nu^2\beta}\cdot \frac{\lambda_2(\beta^{-\rho_2}-1)}{1-\rho_2},\gamma_1^2 \left[1-\frac{\left(R(1,1)-\nu^2R(1,\beta) \right)^2}{1+\nu^2-2\nu^2\beta}\right]\right).$$
  \label{theo1}
 \end{theo}

\noindent \textbf{Remark 1} If we consider again the case $\lambda_1=\lambda_2=0$, the  asymptotic bias is equal to  0 and we can focus on the asymptotic variance.  Consider for convenience the case where  $\beta=1$, then the asymptotic variance becomes $\gamma_1^2[1-(1-\nu^2)R^2(1,1)]$ whereas it is $\gamma_1^2$ for the Hill estimator, meaning that the (relative) variance reduction is equal to $(1-\nu^2)R^2(1,1)$. When, e.g.,  $k_+=2k$, this becomes $\frac{1}{2}R^2(1,1)$. In this case, depending on the value of $R(1,1)\in [0,1]$, the variance reduction can be as large as $50\%$.
 Clearly, in case of tail independence ($R(1,1)=0$), the estimator has the same limiting variance  as the Hill estimator. In such a case a ``better" related variable should be selected.

 \noindent \textbf{Remark 2} Note that in case $\rho_1\neq \rho_2$, we have, since $|A_j|$ is regularly varying at $\infty$ with index $\rho_j$, $j=1,2$, that $\lambda_1=0$ or $\lambda_2=0$. This simplifies the expression for the asymptotic bias. In case $\lambda_2=0$ (which is implied by $\rho_1>\rho_2$) or $\beta=1$ or $\rho_2=0$, the Hill estimator and the adapted estimator have the same asymptotic bias $\lambda_1/(1-\rho_1)$.

\noindent \textbf{Remark 3} It is well-known that choosing a good $k$ is a difficult problem in extreme value theory. We will not address this problem here, but compare for many values of $k$ our adapted estimator and the Hill estimator, see Remark 1, Remark 4, and the simulation section. On the other hand, there  are many methods for choosing the $k$ of the Hill estimator. If one of these methods is adopted, we can choose the same $k$ for our adapted estimator. %Maybe due to the variance reduction we can take $k$  slightly larger.

  \section{Main results: the multivariate case}
  Now we consider a $d$-variate distribution function $F$, with marginals $F_1, \ldots, F_d$ and corresponding tail quantile functions $U_j$, $j=1, \ldots d$; write $F_-$ for the distribution function of the last $d-1$ components of a random vector with distribution function $F$. We assume  that $F$ is in the multivariate max-domain of attraction, that is $F\in D(G) $, with all  extreme value indices $\gamma_1, \ldots, \gamma_d$ positive. Let $R_{ij}$ be the tail copula of the $i$-th and the $j$-th component, $1\leq i,j\leq d, i\neq j$, see~(\ref{tcop}).

  Let $(X_{1},Y_{1,2},\ldots,Y_{1,d}),\ldots,(X_{n},Y_{n,2},\ldots,Y_{n,d}),$ be a $d$-variate random sample from $F$ and let  \\ $(Y_{n+1,2},\ldots,Y_{n+1,d}),\ldots,(Y_{n+m,2},\ldots,Y_{n+m,d})$ be a $(d-1)$-variate random sample from $F_-$, independent of the $d$-variate random sample of size $n$.
%This tail dependence function is estimated as
%$$\hat{R}_{1j}(x,y)=\frac{1}{k}\sum\limits_{l=1}^{n}1_{[X_1\geqslant X_{n-[kx]+1,n}, Y_{1,j}\geqslant Y_{n-[ky]+1,n,j}]}, \quad i=1,\ldots,d \quad j=2,\ldots,d.$$
%The tail dependence function  between $Y_{1,i}$ and $Y_{1,j}$ is
%$$R_{ij}(x,y)=\lim\limits_{t\downarrow 0}\frac{1}{t} \mathbb{P}\Big[1-F_i\big(Y_{1,i} \big) \leq tx,1-F_j\big(Y_{1,j} \big) \leq ty\Big]\quad \quad  \forall  x\geq 0  \text{, and }  y \geq 0.$$
%Hence the estimated version will be
%$$\hat{R}_{ij}(x,y)=\frac{1}{k}\sum\limits_{l=1}^{n}1_{[Y_{l,i}\geqslant Y_{n-[kx]+1,n,i}, Y_{l,j}\geqslant Y_{n-[ky]+1,n,j}]}, \quad i=1,\ldots,d \quad j=2,\ldots,d,$$
%where $Y_{n-[ky]+1,n,i},$ and $Y_{n-[ky]+1,n,j}$ are the $n-[ky]+1$ order statistics from $Y_{l,i}$ and $Y_{l,j}$ respectively $, l=1,2,..,n$.
  Let $\hat{\gamma}_1, \hat{\gamma}_j,$ and $\hat{\gamma}_{j+}$  be the Hill estimators based on $X_{1},\ldots,X_{n},$ $Y_{1,j},\ldots,Y_{n,j},$ and $Y_{1,j},\ldots,Y_{n+m,j},$ $j=2,\ldots,d$, respectively, see (\ref{hil}); here again we replace $k$ with $k_+$ for $\hat{\gamma}_{j+}, j=2, \ldots, d$.
 First we consider the joint asymptotic normality of all the $2d-1$ Hill estimators.
  % When  $k,k_+\rightarrow \infty$ as $n\rightarrow \infty$ assume
%\begin{enumerate}
%\setcounter{enumi}{3}
%    \item[4.] $\lim\limits_{t\to\infty}\frac{\frac{U_j(tx)}{U_j(t)}-x^{\gamma_j}}{A_j(t)}= x^{\gamma_j} \frac{x^{\rho_j-1}}{\rho_j},$ \\ for all $x>0,$ $\rho_j\leqslant0,$ and $A_j$ is a positive or negative function with $\lim \limits_{t\to\infty}A_j(t)=0, j =1,\ldots,d.$
%    \item[5.] $\lim\limits_{n\to\infty} \sqrt{k}A_i(\frac{n}{k})=\lambda_i, \quad \quad \quad \quad i=1,\ldots,d,$ $\quad\quad $ $\lim\limits_{n+m\to\infty} \sqrt{k_+}A_j(\frac{n+m}{k_+})=\lambda_{j+}\quad \quad \quad \quad j=2,\ldots,d.$
%    \end{enumerate}
\begin{prop}
If $F\in D(G)$,  condition (\ref{kaa}) holds, condition (\ref{so}) holds for $j=1, \ldots , d$,  and   $ \sqrt{k}A_j(\frac{n}{k})\to\lambda_j \in \mathbb{R}, j=1,\ldots, d$,
%and $ \sqrt{k_+}A_j(\frac{n+m}{k_+})\to\lambda_{j+} \in \mathbb{R}, j=2,\ldots,d$,
 as $n\to\infty$, then
  \begin{equation}
  \left(\sqrt{k}(\hat{\gamma}_1-\gamma_1),\sqrt{k}(\hat{\gamma}_2-\gamma_2),\sqrt{k_+}(\hat{\gamma}_{2+}-\gamma_2),\ldots,\sqrt{k}(\hat{\gamma}_d-\gamma_d),\sqrt{k_+}(\hat{\gamma}_{d+}-\gamma_{d})\right)\xrightarrow{d} N(\breve\mu_d,\breve\Sigma_d), \label{eq3}
  \end{equation}
where $$ \breve\mu_d=\left(\frac{\lambda_1}{1-\rho_1}, \frac{\lambda_2}{1-\rho_2}, \frac{\lambda_2\beta^{-\rho_2}}{\nu(1-\rho_2)},\ldots,\frac{\lambda_d}{1-\rho_d},  \frac{\lambda_d\beta^{-\rho_d}}{\nu(1-\rho_d)}\right), \quad \quad  $$
$$\breve\Sigma_d=\left[ \begin{array}{ccccccccccc} \gamma_1^2 & R_{12}(1,1) \gamma_1 \gamma_2 & \nu R_{12}(1,\beta) \gamma_1 \gamma_2&.&.&.& R_{1d}(1,1) \gamma_1 \gamma_d & \nu R_{1d}(1,\beta) \gamma_1 \gamma_d \\ R_{12}(1,1) \gamma_1 \gamma_2 &  \gamma_2^2 &\nu\beta \gamma_2^2&.&.&.& R_{2d}(1,1) \gamma_2 \gamma_d&  \nu R_{2d}(1,\beta) \gamma_2 \gamma_d \\. &.&.&.&.&.&.&.&\\ .&.&.&.&.&.&.&.&\\ .&.&.&.&.&.&.&.&\\.&.&.&.&.&.&.&.&\\ R_{1d}(1,1) \gamma_1 \gamma_d  &  R_{2d}(1,1) \gamma_2 \gamma_d &\nu R_{2d}(1,\beta) \gamma_2 \gamma_d &.&.&.& \gamma_d^2&\nu\beta \gamma_d^2\\ \nu R_{1d}(1,\beta) \gamma_1 \gamma_d & \nu R_{2d}(1,\beta) \gamma_2 \gamma_d  & R_{2d}(1,1) \gamma_2 \gamma_d &.&.&.&\nu\beta \gamma_d^2&\gamma_d^2\end{array}\right].$$
\label{prop2}
   \end{prop}
   \begin{cor}
   Under the conditions of Proposition \ref{prop2}, as $n\to\infty$,
   \begin{equation}
   \left(\sqrt{k}(\hat{\gamma}_1-\gamma_1),\sqrt{k}(\hat{\gamma}_{2+}-\hat{\gamma}_2),\ldots,\sqrt{k}(\hat{\gamma}_{d+}-\hat{\gamma}_{d}) \right)\xrightarrow{d}N\Big(\mu_d,\Sigma_d\Big), \label{eq4}
  \end{equation}
   where $\mu_d=\left( \frac{\lambda_1}{1-\rho_1}, \frac{\lambda_2(\beta^{-\rho_2}-1)}{1-\rho_2},\ldots,\frac{\lambda_d(\beta^{-\rho_d}-1)}{1-\rho_d} \right),\, \Sigma_d=\Gamma \Gamma^{T} \circ\,  H$ (``$\,  \circ $" denotes the Hadamard or entrywise product), with $$H= \left[ \begin{array}{ccccccccccc} 1 & h_{12}  &.&.&.& h_{1d} \\ h_{12} &  h &.&.&.& h_{2d} \\ . &.&&&&.\\ .&&.&&&.\\ .&&&.&&.\\.&&&&.&.\\ .&&&&&.\\ h_{1d}& h_{2d} &.&.&.&  h \end{array}\right], \,\, \Gamma=\begin{bmatrix} \gamma_1 \\ \gamma_2 \\.\\.\\.\\ \gamma_d \end{bmatrix},$$  $h=1+\nu^2-2\nu^2\beta,$ $h_{1i}=\nu^2 R_{1i}(1,\beta)-R_{1i}(1,1),$  and $h_{ij}=(1+\nu^2)R_{ij}(1,1)-\nu^2\big(R_{ij}(1,\beta)+R_{ij}(\beta,1)\big)$, $i=2,\ldots,d, $  $j=i+1,\ldots,d.$
  \label{cor2}
   \end{cor}

As in the bivariate case  we   approximate, for $\lambda_j=0, j=1,\ldots,d$, the $d$-variate normal limiting distribution of $(\hat\gamma_1, \hat\gamma_{2+}-\hat\gamma_2, \ldots,   \hat\gamma_{d+}-\hat\gamma_d)$, and estimate the approximated $\frac{1}{k}\Sigma_d$, where for the estimation of $R_{ij}$,  $\hat R_{ij}$ is defined similarly as $\hat R$ in (\ref{rr}). The thus obtained approximated and estimated version of $\frac{1}{k}\Sigma_d$ is denoted by $\frac{1}{k}\hat\Sigma_d$. In this  normal distribution the only unknown parameter is the first component of the mean: $\gamma_1$, cf.\ (\ref{limnor}). Maximizing the corresponding likelihood with respect to $\gamma_1$, we obtain our adapted estimator  for $\gamma_1$:
%\begin{equation}
%f(\tilde{\Gamma}_d,\tilde{\mu}_d,\tilde{\Sigma}_d)=\frac{1}{\sqrt{2\pi^{d+1} |\tilde{\Sigma}_d|}} \exp\left[\frac{-1}{2} (\tilde{\Gamma}_d-\tilde{\mu}_d)^T \tilde{\Sigma}_d^{-1} (\tilde{\Gamma}_d-\tilde{\mu}_d) \right], \label{eq55_}
%\end{equation}  \\  where $$\tilde{\Gamma}_d= \left[\begin{array}{c} \hat{\gamma}_1 \\ \\ \hat{\gamma}_{2+}-\hat{\gamma}_2 \\ \\ \hat{\gamma}_{3+}-\hat{\gamma}_3 \\ .\\.\\  \hat{\gamma}_{d+}-\hat{\gamma}_d \end{array}\right] , \tilde{\mu}_d= \left[\begin{array}{c} \gamma_1 \\ 0 \\ 0 \\ . \\  .\\0 \end{array}\right] , \quad \tilde{\Sigma}_d=\frac{1}{k}\hat{\Gamma} \hat{\Gamma}^{T} \circ \tilde{H}, \quad \tilde{H}=\left[ \begin{array}{ccccccccccc} 1 & \tilde{C}_{12} &.&.&.&.& \tilde{C}_{1d}  \\ \tilde{C}_{12}  &  \tilde{C} &.&.&.&.& \tilde{C}_{2d} \\ .&.&.&.&.&.&.\\ .&&.&&&&.\\ .&&&.&&&.\\.&&&&.&&.\\ .&&&&&.&.\\ \tilde{C}_{1d} & \tilde{C}_{2d}   &.&.&.&.&  \tilde{C}\end{array}\right],\hat{\Gamma}=\begin{bmatrix} \hat{\gamma}_1 \\ \hat{\gamma}_{2+} \\.\\.\\.\\\hat{\gamma}_{d+} \end{bmatrix},$$ $ \tilde{C}=1+\frac{k}{k_+}-\frac{2n}{n+m},$ $\tilde{C}_{1i}=\frac{k}{k_+} \hat{R}_{1i}(1,\tilde{\beta})-\hat{R}_{1j}(1,1),$  and $\tilde{C}_{ij}=\hat{R}_{ij}(1,1)+\frac{k}{k_+}\big(\hat{R}_{ij}(1,1)-2\hat{R}_{ij}(1,\tilde{\beta})\big)$, $i=2,...,d,$  $j=i+1,..d$. \noindent We can apply the MLE on \ref{eq55_} to obtain estimator for $\gamma_1$. The MLE yields to
$$\hat{\gamma}_{1,d}=\hat{\gamma}_1+\sum\limits_{j=2}^{d}\frac{\hat{\Sigma}_{1j}^{-1}}{\hat{\Sigma}_{11}^{-1}}(\hat{\gamma}_{j+}-\hat{\gamma}{_j}),$$
%+\frac{\tilde{\Sigma}_{13}^{-1}}{\tilde{\Sigma}_{11}^{-1}}(\hat{\gamma}_{3+}-\hat{\gamma}{_3})+\ldots+\frac{\tilde{\Sigma}_{1d}^{-1}}{\tilde{\Sigma}_{11}^{-1}}(\hat{\gamma}_{d+}-\hat{\gamma}{_d}),$$ \\
 where $A^{-1}_{ij}$ denotes the entry in the $i^{th}$ row and $j^{th}$ column of the inverse of the matrix $A$. Using, in the obvious notation, $\hat\Sigma_d=\hat\Gamma \hat \Gamma^{T} \circ \, \hat H$ (see above),
  %According to the Hadamard product $ \tilde{\Sigma}^{-1}=\left( \frac{1}{k} \hat{\Gamma} \hat{\Gamma}^{T}\right)^{\circ(-1)} \circ \tilde{H}^{-1}$, where $\left(\frac{1}{k} \hat{\Gamma} \hat{\Gamma}^{T}\right)^{\circ(-1)}$ is the Hadamard inverse (i.e. the element wise inverse) of $\left(\frac{1}{k} \hat{\Gamma} \hat{\Gamma}^{T}\right)$. Hence  $\tilde{\Sigma}^{-1}_{11}=k \frac{\tilde{H}^{-1}_{11}}{\hat{\gamma}_{1}^2},$ and $\tilde{\Sigma}^{-1}_{1i}=k \frac{\tilde{H}^{-1}_{1i}}{\hat{\gamma}_{1} \hat{\gamma}_{i+}}, i=2,\ldots,d,$
   we can  rewrite our adapted estimator as
\begin{equation}\label{ahm} \hat {\gamma}_{1,d}
 %\hat{\gamma}_1+\frac{\hat{\gamma}_{1}}{\hat{\gamma}_{2+}} \frac{\tilde{H}_{12}^{-1}}{\tilde{H}_{11}^{-1}}(\hat{\gamma}_{2+}-\hat{\gamma}{_2})+\frac{\hat{\gamma}_{1}}{\hat{\gamma}_{3+}} \frac{\tilde{H}_{13}^{-1}}{\tilde{H}_{11}^{-1}}(\hat{\gamma}_{3+}-\hat{\gamma}{_3})+\ldots+\frac{\hat{\gamma}_{1}}{\hat{\gamma}_{d+}} \frac{\tilde{H}_{1d}^{-1}}{\tilde{H}_{11}^{-1}}(\hat{\gamma}_{d+}-\hat{\gamma}{_d})
 =\hat{\gamma}_1+  \sum_{j=2}^{d} \frac{\hat{\gamma}_1}{\hat{\gamma}_{j+}} \frac{\hat{H}_{1j}^{-1}}{\hat {H}_{11}^{-1}} (\hat{\gamma}_{j+}-\hat{\gamma}{_j}).\end{equation}
% where $\hat{H}^{-1}_{1j}$ is the entry in the $1^{st}$ row and $j^{th}$ column of the matrix $\hat{H}^{-1}$.

 \begin{theo}
Assume $H$ is invertible. Then under the conditions of Proposition \ref{prop2}, as $n\to\infty$,
 \begin{equation}
  \sqrt{k}(\hat{\gamma}_{1,d}-\gamma_1)\xrightarrow{d}N\left(\frac{\lambda_1}{1-\rho_1}+ \sum \limits_{j=2}^{d}\frac{\gamma_1}{\gamma_j}\frac{H_{1j}^{-1}}{H_{11}^{-1}}\frac{\lambda_j(\beta^{-\rho_j}-1)}{1-\rho_j},\sigma^2\right), \label{eq5_}
  \end{equation}
 where
 %$$ABias(\tilde{\gamma}_{1,d})=\frac{\lambda_1}{1-\rho_1}+ \sum \limits_{j=2}^{d}\frac{\gamma_1}{\gamma_j}\frac{H_{1j}^{-1}}{H_{11}^{-1}}\frac{\nu\lambda_{j+}-\lambda_j}{1-\rho_j},$$
 \begin{multline*}
\sigma^2=\gamma_1^2 \Big(1-\frac{1}{(H_{11}^{-1})^2}\Big[2 H_{11}^{-1}\sum\limits_{j=2}^{d}H_{1j}^{-1}[R_{1j}(1,1)-\nu^2R_{1j}(1,\beta)]-[1+\nu^2-2\nu^2\beta]\sum\limits_{j=2}^{d}
(H_{1j}^{-1})^2\\-2\sum\limits_{i=2}^{d}\sum\limits_{j>i}^{d}[(1+\nu^2)R_{ij}(1,1)-\nu^2\big(R_{ij}(1,\beta)+R_{ij}(\beta,1) \big)\big]H_{1i}^{-1} H_{1j}^{-1}\Big]\Big).
 \end{multline*}
% where $H_{1i}^{-1}$ is the entry in the $1^{st}$ row and $i^{th}$ column of the matrix $H^{-1}$.
\label{theo2}
 \end{theo}

\noindent \textbf{Remark 4} We have seen that in the bivariate case for  $\beta=1$ and $\nu^2=\frac{1}{2}$ the  reduction in asymptotic variance is equal to $\frac{1}{2}R^2(1,1)$. For, e.g., $R(1,1)=0.8$, this becomes $0.320$.
Now consider the trivariate case with the same values for $\beta$ and $\nu^2$ and with (also) $R_{12}(1,1)=R_{13}(1,1)=0.8$, but $R_{23}(1,1)=0.4$. Then the reduction in asymptotic variance, see the next section, becomes much larger: $0.457$. In other words, adding a third variable that has the same (as the second variable) tail copula value at (1,1)  with the variable of interest and does not have a high tail dependence with the second variable  reduces the asymptotic variance much more than when using only one related variable.

\section{Simulation study}

In this section we will perform a simulation study in order to compare the finite sample behavior of the adapted  estimator and the Hill estimator. We will consider 6 bivariate distributions and 8 trivariate distributions and 3 different pairs $(n,m)$. Every setting is replicated 10,000 times.

To be precise, we consider the Cauchy distribution restricted to the first quadrant/octant in dimensions $d=2$ and $d=3$. This Cauchy density is proportional to
$$(1+xS^{-1}x^T)^{-(1+d)/2},$$ where  the $2\times2$ or $3\times 3$ scale matrix $S$ has 1 as diagonal elements and $s$ as off-diagonal elements, but when $d=3$ we take $S_{23}=S_{32}=r$. For $s$ we take the values 0, 0.5, and 0.8, respectively. When $d=3$ we take $r=s$, but for $s=0.5$  and $s=0.8$ we also take $r=0$  and $r=0.3$, respectively.
 We will also consider the bi- and trivariate logistic distribution function  with standard Fr\'echet marginals:
 $$F(x_1,\ldots, x_d)=\exp\left\{ -\left(x_1^{-1/\theta}+ \ldots +x_d^{-1/\theta}\right)^\theta\right\}, \quad x_1>0, \ldots, x_d>0; \quad d=2 \mbox{ or } d=3.$$
 For $\theta$ we take the values 0.1, 0.3, and 0.5, respectively.

 We use the following values for $n$, $m$, and $k$: \\
  $\bullet$ $n=1000$, $m=500$, and $k=100$,\\
  $\bullet$ $n=1000$, $m=1000$, and $k=100$, \\
 $\bullet$ $n=500$, $m=1000$, and $k=50$.\\
Then we choose $k_+$ according to
\begin{equation}
\frac{k}{k_+}=\frac{n}{n+m}. \label{eq15}
\end{equation}
%which yields $1+\frac{k}{k_+}-2\frac{n}{n+m}=1-\frac{n}{n+m}$ and  $\beta=1$.

 In case $d=2$, using (\ref{eq15}), our adapted estimator in (\ref{eq2}) specializes to
 $$\hat{\gamma}_{1,2}=\hat{\gamma}_1+\frac{\hat{\gamma}_1}{\hat{\gamma}_{2+}}\hat{R}(1,1)(\hat{\gamma}_{2+}-\hat{\gamma}{_2}),$$
and the asymptotic variance in Theorem \ref{theo1} becomes
$\gamma_1^2 \Big(1-(1-\nu^2)R^2(1,1)\Big)$.
%%where $0 \leq (1-\nu^2)R^2(1,1) \leq 1.$ \\  \\
When $d=3$,  using (\ref{eq15}), our adapted estimator in (\ref{ahm})  can be rewritten as
$$\hat{\gamma}_{1,3}=\hat{\gamma}_1+\frac{\hat{\gamma}_1}{\hat{\gamma}_{2+}}\frac{\hat{R}_{12}(1,1)-\hat{R}_{13}(1,1)\hat{R}_{23}(1,1)}{1-\hat{R}^2_{23}(1,1)}
(\hat{\gamma}_{2+}-\hat{\gamma}{_2})+\frac{\hat{\gamma}_1}{\hat{\gamma}_{3+}}\frac{\hat{R}_{13}(1,1)-\hat{R}_{12}(1,1)\hat{R}_{23}(1,1)}{1-\hat{R}^2_{23}(1,1)}(\hat{\gamma}_{3+}-\hat{\gamma}{_3}),$$
and the  asymptotic variance in Theorem \ref{theo2} specializes to
$$\sigma^2=\gamma_1^2\left(1-(1-\nu^2) \left(\frac{R^2_{12}(1,1)+R^2_{13}(1,1)-2R_{12}(1,1)R_{13}(1,1)R_{23}(1,1)}{1-R^2_{23}(1,1)}\right)\right).$$
%where $0 \leq \frac{R^2_{12}(1,1)+R^2_{13}(1,1)-2R_{12}(1,1)R_{13}(1,1)R_{23}(1,1)}{1-R^2_{23}(1,1)} \leq 1$.

Tables \ref{tab} and \ref{tab2} show the  (empirical percentages of) variance reduction as discussed below Theorems \ref{prop1} and \ref{prop2}, and above. We see that the variance reduction ranges from about 10\%
\begin{table}[h!]
\centering
\vspace{ .5 cm}
\begin{tabular}{lcccccccc} \hlinewd{.64pt}
\multicolumn{1}{c}{\text{}} & \multicolumn{3}{c}{$d=2$} & \multicolumn{5}{c}{$d=3$} \\
\hline
%\cmidrule(r){4-5}
\text{ } & $s=0$ & $s=0.5$   & $s=0.8$  & $s=0$ & $s=0.5$   & $s=0.5$&$s=0.8$&$s=0.8$ \\
\text{ } &  &    &   & $r=0$ & $r=0.5$ &$r=0$   & $r=0.8$&$r=0.3$ \\
\hline
$n=1000, m=500$ & 10.5\% & 12.1\% & 16.7\% &12.7\%  & 17.6\%& 19.0\%&21.9\%& 25.8\%   \\
$n=1000, m=1000$ & 16.1\% &20.9\% &27.2\% & 19.8\%  & 26.0\%&30.1\% &32.2\%& 39.3\%   \\
$n=500, m=1000$ & 20.8\% & 28.3\% & 37.3\% &26.6\%& 34.3\% &37.0\%&42.6\%&   52.7\%   \\
\hlinewd{.72pt}
\end{tabular}
  \caption{Empirical variance reduction for the Cauchy distribution}
  \label{tab}
\end{table}
\begin{table}[h!]
\centering
\begin{tabular}{lcccccc}
 \hlinewd{.64pt}
%\toprule
\multicolumn{1}{c}{\text{}} & \multicolumn{3}{c}{$d=2$} & \multicolumn{3}{c}{$d=3$} \\
\hline
%\cmidrule(r){2-3}
%\cmidrule(r){4-5}
\text{ } & $\theta=0.1$ & $\theta=0.3$   & $\theta=0.5$  & $\theta=0.1$ & $\theta=0.3$   & $\theta=0.5$ \\ \hline
%\midrule
$n=1000, m=500$ & 26.6\% & 18.1\% & 9.1\% &27.4\%  & 20.4\%&13.2\% \\
$n=1000, m=1000$ & 41.7\% &27.7\% & 15.2\% & 44.5\%  & 33.1\% & 20.6\%\\
$n=500, m=1000$ & 55.6\% & 36.3\% & 21.7\% & 57.0\%& 42.1\% &26.1\% \\
\hlinewd{.72pt}
%\bottomrule
\end{tabular}
  \caption{Empirical variance reduction for the logistic  distribution}
  \label{tab2}
\end{table}
to more than 50\%, that is, our adapted estimator  yields much better results than the Hill estimator. A stronger  tail dependence between the variable of interest and the related variable(s) yields a larger variance reduction. In case $d=3$, due to the exchangeability of the components of the logistic distribution, a stronger tail dependence between the variable of interest and the related variables, yields  also  a stronger tail dependence between the two related variables and hence increasing the dimension from 2 to 3 does not help that much, but in case of the Cauchy distribution with $r<s$ we see a large improvement when adding the third variable. Comparing the numbers in the table with the (not presented) theoretical asymptotic reductions shows that the empirical numbers are about the same but slightly smaller, partly due to the variability of the  tail copula estimators, which does not show up in the asymptotic variance.

\begin{figure}[h!]
\begin{center}\label{bpah}
\includegraphics[width=.8\textwidth, height=6cm]{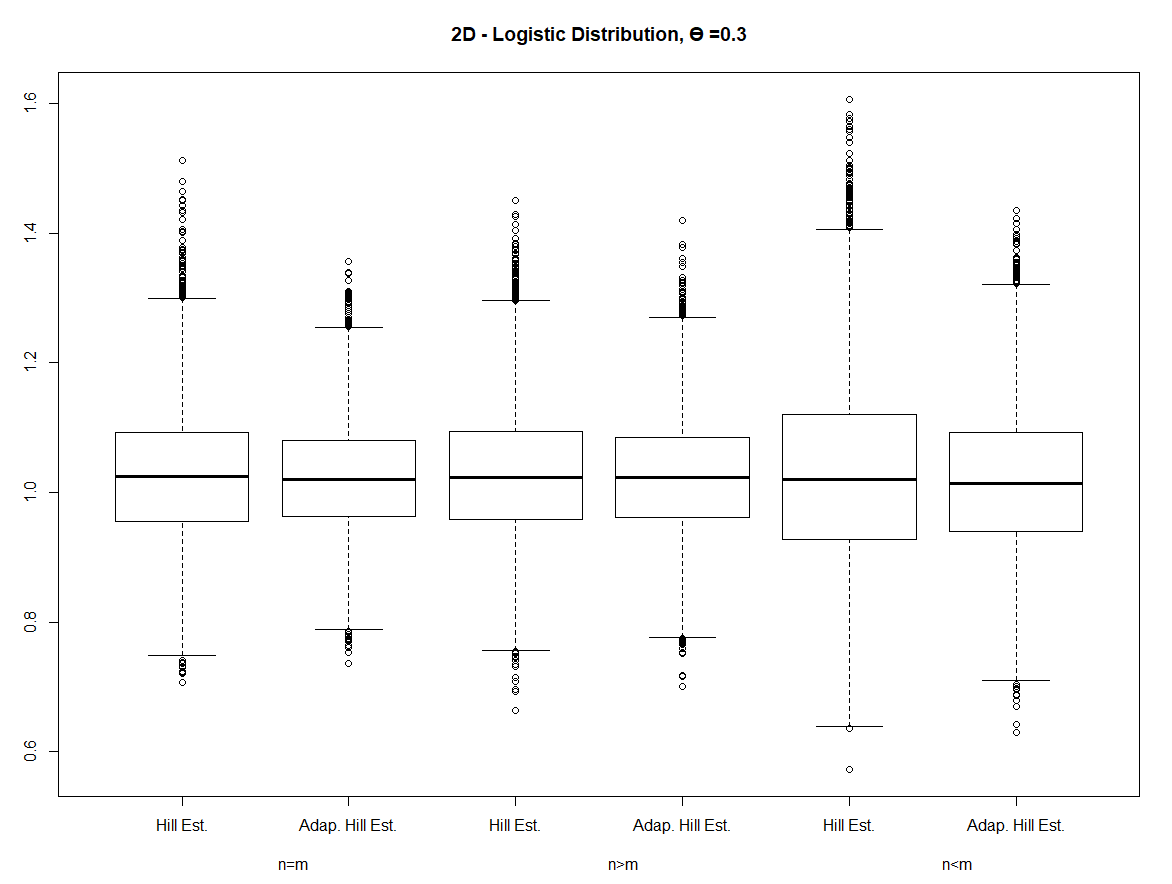}\\ \vspace{0.5 cm}
\includegraphics[width=.8\textwidth, height=6cm]{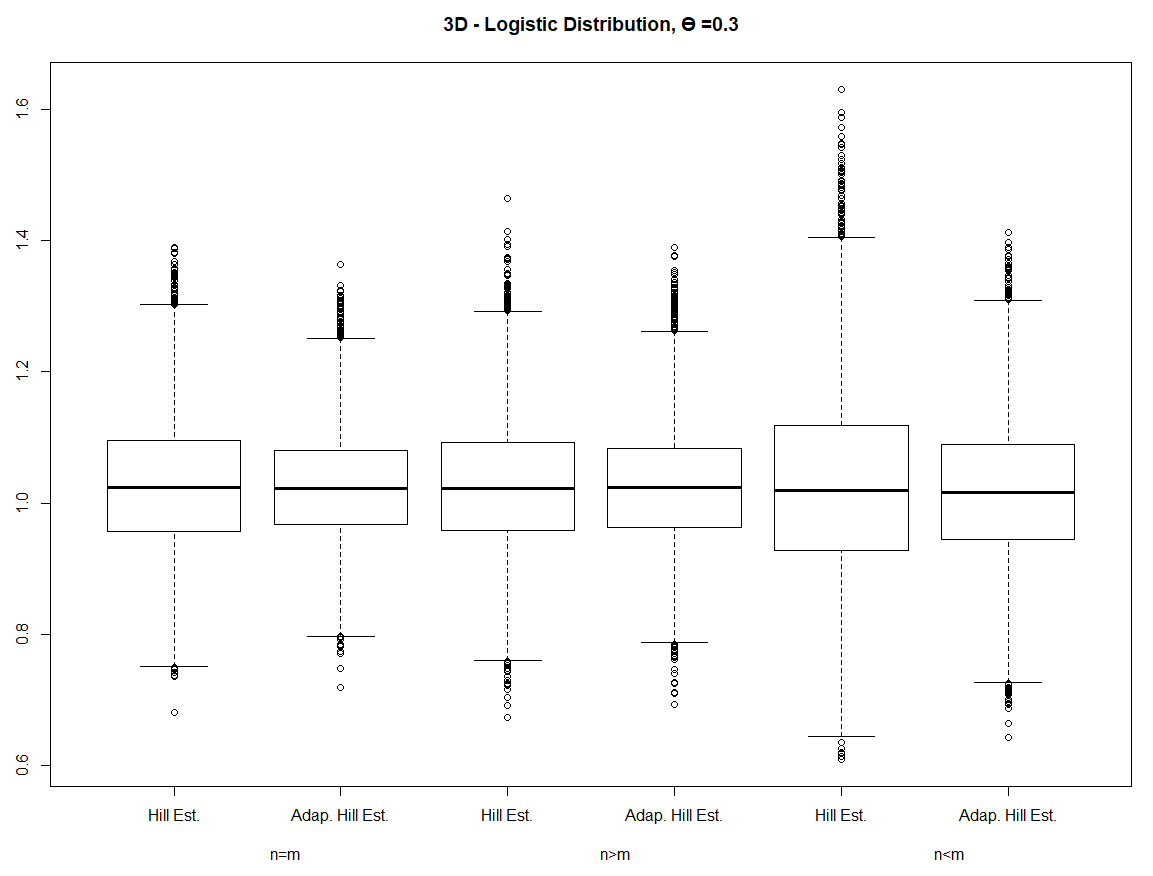}
\caption{Boxplots of 10,000 estimates of $\gamma_1$ based on the logistic distribution with $\theta=0.3$.}
\end{center}
\end{figure}

In order to briefly investigate not only the variance of the estimators,  but their full behavior  we also present boxplots for our estimator and the Hill estimator corresponding to, as an example,  the logistic distribution with $\theta=0.3$ in dimensions 2 and 3, with the $(n,m)$-settings as before (Figure 1). Again we see that our adapted estimator outperforms the Hill estimator: the boxes are smaller and there are less outlying estimates.

\section{Application}

We apply our bivariate adapted estimator to  financial losses (in US\$) due  to earthquakes; the related variable is the corresponding  energy released. The aim of this application is to assess the tail heaviness of the loss distribution and also to estimate a very high quantile of the losses.

The earthquakes concern  29 countries\footnote{Algeria, Burma, Chile, China, Ecuador, El Salvador, Germany, Greece, Haiti, Iceland, India, Indonesia, Iran, Italy, Japan, Mexico, Morocco, Nepal, New Zealand, Nicaragua, Pakistan, Philippines, Russia, Taiwan, Tajikistan, Tanzania, Thailand, Turkey, United States.}.
The data are provided by  the National Oceanic and Atmospheric Administration (NOAA).
Ignoring tsunami  losses, we consider  the financial losses of categories at least  ``moderate" for  the       time period from 1993 through 2017.
 We used linear regression analysis per country   for imputation of  missing loss values, with ``number of deaths due to the earthquake" and  ``severity of the financial loss" (a categorical variable)  as independent variables.
 We also corrected the financial losses for inflation. The highest loss in the data set is  US\$~36$\times 10^9$.
We obtained the related Richter scale magnitude $M$ of the earthquakes  for the much longer period 1940 through 2017. The energy $E$ released by earthquakes (in megajoules) is given by
$E=2\times 10^{1.5(M-1)}$; \cite{lay1995modern}.

 \begin{figure}[h!]
\begin{center}\label{AHP}\vspace{0 cm}
 \includegraphics[width=10 cm, height= 5.5 cm]{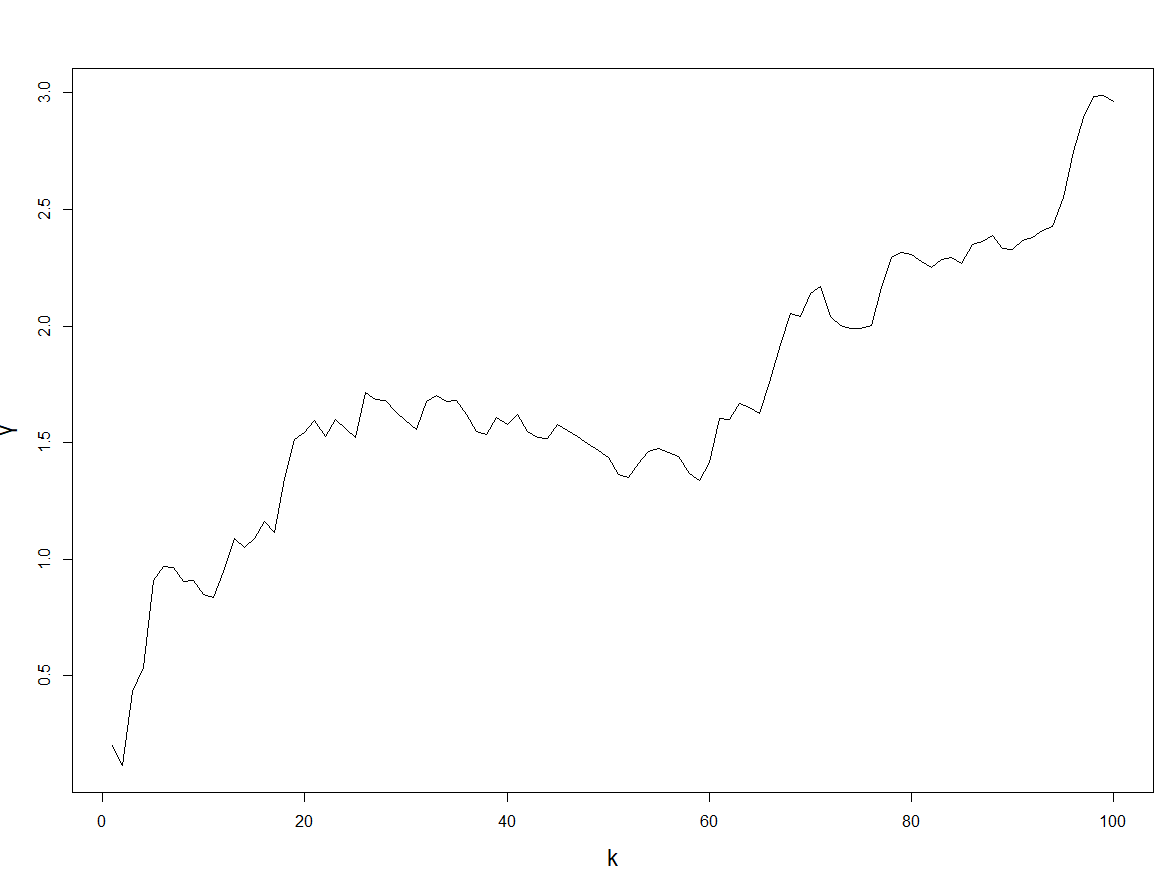}
 \caption{Adapted Hill estimator  of the  financial losses of the earthquakes}
\end{center}
\end{figure}

We have  $n=330$ and $m=512$.      Figure 2 shows a plot of the adapted Hill estimator against $k$, with $k_+$ based on (\ref{eq15}).
      We  take the average value of the estimates  over the region $k=40, \ldots, 60 $. This yields the average Hill estimate  $\hat{\gamma}_1=1.504 $
and our final average estimate of $\gamma_1$, which is somewhat lower than the Hill estimate:
 $$ \hat{\gamma}_{1,2}=1.465.$$
 Both estimates indicate that the loss distribution has a very heavy right tail.

 We also estimate the high quantile $F_1^{-1}(1-p)$ of the loss distribution for $p=\frac{1}{n}=\frac{1}{330}$.  %using  both estimates of the extreme value index, respectively.
 This high quantile is estimated as  usual (see, e.g., page 138 of \cite{de2006extreme}) with $$X_{n-k,n}\left(\frac{k}{np}\right)^{\hat{\gamma}},$$
 where ${\hat{\gamma}}$ is the Hill estimator or the adapted Hill estimator (and $k=40, \ldots, 60$).
 %We have  $X_{n-k,n}=$ \\ US\$ 279$\times 10^6$.
This yields for the average high quantile estimate US\$ 130$\times 10^9$ when we use the Hill estimates and US\$~113$\times 10^9$ when we use our estimates of $\gamma_1$, which is a reduction of 17 billion dollars. This shows that, from an insurer's perspective, improved (that is, less variable) estimation of the extreme value index can lead to huge changes in high quantiles, here the 25 year return level. The lower estimate we obtain indicates less risk for (re)insurers.

\section{Proofs}
\textbf{Proof of Proposition \ref{prop2}}
Let $C$ be a copula corresponding to the distribution function of $(-X_1,-Y_{1,2},\ldots,-Y_{1,d})$ and let $C_-$ be the distribution function of the last $d-1$ components of a random vector with distribution function $C$.  Let $(V_{1,1},V_{1,2},\ldots,V_{1,d}),\ldots,(V_{n,1},V_{n,2},\ldots,V_{n,d})$ be a random sample of size $n$ from $C$ and let $(V_{n+1, 2},\ldots,V_{n+1, d}),\ldots,(V_{n+m,2},\ldots,V_{n+m,d})$ be a random sample of size $m$ from $C_-$,  independent of the random  sample from $C$. Clearly all the $V_{i,j}$ have a uniform-(0,1) distribution.
Write $X_i=F_1^{-1}(1-V_{i,1}), \,   i=1,\ldots,n$,  and $Y_{l,j}=F_j^{-1}(1-V_{l,j}),\,  l=1,\ldots,n+m,\, j=2,\ldots,d.$ Then  $(X_{1},Y_{1,2},\ldots,Y_{1,d}),\ldots,(X_{n},Y_{n,2},\ldots,Y_{n,d}),$ and $(Y_{n+1,2},\ldots,Y_{n+1,d}),\ldots,(Y_{n+m,2},\ldots,Y_{n+m,d})$ have the distributions as specified in the beginning of Section 3.

Consider the univariate empirical distribution functions  $\Gamma_{n,j}(s)=\frac{1}{n}\sum_{i=1}^n 1_{[0,s]}(V_{i,j}),\, \\ 0\leq s\leq1,\,  j=1,2,\ldots,d,$ and
$\Gamma_{n+m,j}(t)=\frac{1}{n+m}\sum_{l=1}^{n+m} 1_{[0,t]}(V_{l,j}),\,  0\leq t\leq 1,\,   j=2,\ldots,d,$ and the corresponding uniform tail empirical processes
\begin{eqnarray*}
&& w_{n,j}(s)=\frac{n}{\sqrt{k}}\left[\Gamma_{n,j}\left(\frac{k}{n}s\right)-\frac{k}{n}s\right], \quad 0\leq s \leq 1, \\
&& w_{n+m,j}(t)=\frac{n+m}{\sqrt{k_+}}\left[\Gamma_{n+m,j}\left(\frac{k_+}{n+m}t\right)-\frac{k+}{n+m}t\right], \quad\ 0\leq t \leq 1.
\end{eqnarray*}
Now define the Gaussian vector of  processes $\left(W_1,\dots, W_{2d-1}\right)$, where $W_j,\, j=1, \ldots, 2d-1,$ is a standard Wiener process on $[0,1]$,
and the covariances are as follows:
\begin{eqnarray}
&& Cov(W_i(s),W_j(t))= R_{ij}(s,t), \quad 0 \leq s,t \leq 1,  \quad  1 \leq i < j \leq d,\nonumber\\
&& Cov(W_i(s),W_j(t))= \nu R_{i,j-d+1}(s,\beta t), \quad 0 \leq s,t \leq 1,  \quad  1 \leq i  \leq d, \,  d+1 \leq j \leq 2d-1, j\neq i+d-1,\nonumber\\
&& Cov(W_i(s),W_{i+d-1}(t))= \nu (s\wedge \beta t), \quad 0 \leq s,t \leq 1,\, 2\leq  i \leq d.\nonumber\\
&& Cov(W_i(s),W_j(t))= R_{i-d+1,j-d+1}(s,t), \quad 0 \leq s,t \leq 1,  \quad  d+1 \leq i <j \leq 2d-1. \label{cova}
\end{eqnarray}

Let $I$ denote the identity function on $[0,1]$. Then we have   on $(D[0,1])^{2d-1}$, for $0\leq \delta < \frac{1}{2}$, as $n \rightarrow \infty$,
\begin{equation}
\left(\frac{w_{n,1}}{I^\delta},\ldots,\frac{w_{n,d}}{I^\delta},\frac{w_{n+m,2}}{I^\delta},\ldots,\frac{w_{n+m,d}}{I^\delta}\right)\xrightarrow{d} \left(\frac{W_1}{I^\delta},\ldots,\frac{W_d}{I^\delta},\frac{W_{d+1}}{I^\delta},\ldots,\frac{W_{2d-1}}{I^\delta}\right).
\label{eq6}
\end{equation}
For the proof of this statement, note that the convergence and tightness of every component is well-known, see Corollary 4.2.1 in  \cite{cs} or Theorem 3 in \cite{einm}. This also yields  the tightness of the entire vector on the left-hand side. It remains to prove the convergence of the finite-dimensional distributions (without the $I^\delta$), which follows from the  (general) multivariate central limit theorem. It suffices to compute the limits of the covariances: we perform this computation for the second formula in (\ref{cova}); the other three formulas there are essentially special cases of that one. We have
%\begin{multline*}
%Cov(V_{n,1}(r),V_{n,j}(s))=Cov(\frac{1}{\sqrt{k}}\sum\limits_{i=1}^n 1_{[0,\frac{k}{n}r]}(V_{11}),\frac{1}{\sqrt{k}}\sum\limits_{i=1}^n 1_{[0,\frac{k}{n}s]}(V_{1j}))\\=\frac{n}{k}Cov( 1_{[0,\frac{k}{n}r]}(V_{11}), 1_{[0,\frac{k}{n}s]}(V_{1j}))=\frac{n}{k}\big[P(V_{11} \leqslant \frac{k}{n}r,V_{1j} \leqslant \frac{k}{n}s)- \frac{k^2}{n^2}rs\big]= \frac{n}{k}P(V_{11} \leqslant \frac{k}{n}r,V_{1j} \leqslant \frac{k}{n}s)- \frac{k}{n}rs\rightarrow R_{1j}(r,s)\\=Cov(W_1(r),W_j(s))\quad j=2,..,d.
%\end{multline*}
%\begin{multline*}
%Cov(V_{n,l}(r),V_{n,j}(s))=Cov(\frac{1}{\sqrt{k}}\sum\limits_{i=1}^n 1_{[0,\frac{k}{n}r]}(V_{il}),\frac{1}{\sqrt{k}}\sum\limits_{i=1}^n 1_{[0,\frac{k}{n}s]}(V_{ij}))\\=\frac{n}{k}Cov( 1_{[0,\frac{k}{n}r]}(V_{1l}), 1_{[0,\frac{k}{n}s]}(V_{1j}))=\frac{n}{k}\big[P(V_{1l} \leqslant \frac{k}{n}r,V_{1j} \leqslant \frac{k}{n}s)- \frac{k^2}{n^2}rs\big]= \frac{n}{k}P(V_{1l} \leqslant \frac{k}{n}r,V_{1j} \leqslant \frac{k}{n}s)- \frac{k}{n}rs\rightarrow R_{lj}(r,s)\\=Cov(W_l(r),W_j(s))\quad l,j=2,..,d \quad \text{and} \quad l \neq j.
%\end{multline*}
\begin{eqnarray*}
&&Cov(w_{n,i}(s),w_{n+m,j-d+1}(t))=Cov\left(\frac{1}{\sqrt{k}}\sum\limits_{l=1}^n 1_{[0,\frac{k}{n}s]}(V_{l,i}),\frac{1}{\sqrt{k_+}}\sum\limits_{l=1}^{n+m} 1_{[0,\frac{k_+}{n+m}t]}(V_{l,j-d+1})\right)\\
&&=Cov\left(\frac{1}{\sqrt{k}}\sum\limits_{l=1}^n 1_{[0,\frac{k}{n}s]}(V_{l,i}),\frac{1}{\sqrt{k_+}}\sum\limits_{l=1}^{n} 1_{[0,\frac{k_+}{n+m}t]}(V_{l,j-d+1})\right)\\
&&=\frac{n}{\sqrt{kk_+}}Cov( 1_{[0,\frac{k}{n}s]}(V_{1,i}), 1_{[0,\frac{k_+}{n+m}t]}(V_{1,j-d+1}))\\
&&=  \frac{n}{\sqrt{kk_+}}\left[\mathbb{P}\left(V_{1,i} \leq \frac{k}{n}s,V_{1,j-d+1} \leq \frac{k_+}{n+m}t\right)- \frac{kk_+}{n(n+m)}st\right]\\
&&=\sqrt{\frac{k}{k_+}}\left[\frac{n}{k}\mathbb{P}\left(V_{1,i} \leq \frac{k}{n}s,V_{1,j-d+1} \leq \frac{k}{n} \frac{n}{k} \frac{k_+}{n+m}t\right)- \frac{k_+}{n+m}st\right]\\
&&\quad \quad  \rightarrow \nu R_{i,j-d+1}(s,\beta t)=Cov(W_i(s),W_{j}(t)).
\end{eqnarray*}
Hence (\ref{eq6}) is established.
%\begin{multline*}
%Cov(V_{n,j}(r),V_{n+m,j}(s))=Cov(\frac{1}{\sqrt{k}}\sum\limits_{i=1}^n 1_{[0,\frac{k}{n}r]}(V_{ij}),\frac{1}{\sqrt{k_+}}\sum\limits_{i=1}^{n+m} 1_{[0,\frac{k_+}{n+m}s]}(V_{ij}))\\=Cov(\frac{1}{\sqrt{k}}\sum\limits_{i=1}^n 1_{[0,\frac{k}{n}r]}(V_{ij}),\frac{1}{\sqrt{k_+}}\sum\limits_{i=1}^{n} 1_{[0,\frac{k_+}{n+m}s]}(V_{ij})) =\frac{n}{\sqrt{kk_+}}Cov( 1_{[0,\frac{k}{n}r]}(V_{1j}), 1_{[0,\frac{k_+}{n+m}s]}(V_{1j}))\\= \frac{n}{\sqrt{kk_+}}\big[P(V_{1j} \leqslant \frac{k}{n}r,V_{1j} \leqslant \frac{k_+}{n+m}s)- \frac{kk_+}{n(n+m)}rs\big]=\sqrt{\frac{k}{k_+}}\big[\frac{n}{k}P(V_{1j} \leqslant \frac{k}{n}r,V_{1j} \leqslant \frac{k}{n} \frac{n}{k} \frac{k_+}{n+m}s)- \frac{k_+}{n+m}rs\big] \\ \rightarrow \nu (r \wedge \beta s)= Cov(W_j(r),W_{j+d-1}(s)) \quad j=2,\ldots,d.
%\end{multline*}
%\begin{multline*}
%Cov(V_{n+m,j}(r),V_{n+m,l}(s))=Cov(\frac{1}{\sqrt{k_+}}\sum\limits_{i=1}^{n+m} 1_{[0,\frac{k_+}{n+m}r]}(V_{ij}),\frac{1}{\sqrt{k_+}}\sum\limits_{i=1}^{n+m} 1_{[0,\frac{k_+}{n+m}s]}(V_{il}))\\ =\frac{n+m}{k_+}Cov( %1_{[0,\frac{k_+}{n+m}r]}(V_{1j}), 1_{[0,\frac{k_+}{n+m}s]}(V_{1l}))= \frac{n+m}{k_+}\big[P(V_{1j} \leqslant \frac{k_+}{n+m}r,V_{1l} \leqslant \frac{k_+}{n+m}s)- \frac{k_+^2}{(n+m)^2}rs\big]\\=\frac{n+m}{k_+}P(V_{1j} %\leqslant \frac{k_+}{n+m}r,V_{1l} \leqslant \frac{k+}{n+m} s)- \frac{k_+}{n+m}rs \\ \rightarrow R_{jl}(r,s)= Cov(W_j(r),W_{j+d-1}(s)) \quad j,l=2,\ldots,d \quad j\neq l.
%\end{multline*}
%\\

According to   \cite{de2006extreme}, Chapter 5 and Theorem 2.3.9,  we have, as $n\to \infty$,
\begin{equation*}
\sqrt{k}(\hat{\gamma}_j-\gamma_j)= -\gamma_j \left(w_{n,j}(1)-\int_{0}^{1}\frac{w_{n,j}(u)}{u} du\right)+\frac{\lambda_j}{1-\rho_j}+o_p(1),\quad  j=1, \ldots, d.
\end{equation*}
Using that $|A_j|$ is regularly varying at $\infty$ with index $\rho_j$, we get similarly
\begin{equation*}
\sqrt{k}(\hat{\gamma}_{j+}-\gamma_j)= -\gamma_j \left(w_{n+m,j}(1)-\int_{0}^{1}\frac{w_{n+m,j}(u)}{u} du\right)+\frac{\lambda_j(\beta^{-\rho_j}-1)}{1-\rho_j}+o_p(1),\quad  j=2, \ldots, d.
\end{equation*}
Combining all these with  (\ref{eq6}) we obtain
\begin{eqnarray*}
&& \left(\sqrt{k}(\hat{\gamma}_1-\gamma_1),
\ldots,\sqrt{k}(\hat{\gamma}_d-\gamma_d),\sqrt{k_+}(\hat{\gamma}_{2+}-\gamma_2),\ldots,\sqrt{k_+}(\hat{\gamma}_{d+}-\gamma_d)\right) \\
&&\xrightarrow{d}\left( -\gamma_1 \left(W_1(1)-\int\limits_{0}^{1}\frac{W_1(u)}{u} du\right)+\frac{\lambda_1}{1-\rho_1},
\ldots,
-\gamma_d \left(W_{d}(1)-\int\limits_{0}^{1} \frac{W_d(u)}{u} du  \right)
+\frac{\lambda_d}{1-\rho_{d}},\right.\\
&&\hspace{-1 cm}\left.-\gamma_{2} \left(W_{d+1}(1)-\int\limits_{0}^{1} \frac{W_{d+1}(u)}{u} du\right)+\frac{\lambda_2(\beta^{-\rho_2}-1)}{1-\rho_{2}},
\ldots,
 -\gamma_{d} \left(W_{2d-1}(1)-\int\limits_{0}^{1}\frac{W_{2d-1}(u)}{u} du\right)+\frac{\lambda_d(\beta^{-\rho_d}-1)}{1-\rho_{d}} \right).
\end{eqnarray*}
It is immediate and well-known that this yields the mean vector and the variances as in the proposition. (Note that the components of the left-hand side there are listed in a different order.) It remains to derive the covariances. Again we only consider the case where $1 \leq i  \leq d,$ $  d+1 \leq j \leq 2d-1, j\neq i+d-1$. The other cases are easier and essentially special cases of this one. We have
\begin{eqnarray}\label{co}
&&\!\!\!\!\!\!\!\!\!\!\!\!\!\! Cov\left(-\gamma_i \left(W_i(1)-\int\limits_{0}^{1} \frac{W_i(u)}{u} du\right)+\frac{\lambda_i}{1-\rho_i},-\gamma_{j-d+1} \left(W_j(1)-\int\limits_{0}^{1}  \frac{W_j(v)}{v} dv\right)+\frac{\lambda_{j-d+1}(\beta^{-\rho_{j-d+1}}-1)}{1-\rho_{j-d+1}}\right)\nonumber \\ \nonumber
&&\!\!\!\!\!\!\!\!\!\!\!\!\!\! =  \gamma_i \gamma_{j-d+1} \big[E(W_i(1)W_j(1))+
  \int\limits_{0}^{1} \int\limits_{0}^{1}\frac{E \big(W_i(u)W_j(v) \big)}{uv} dudv-\int\limits_{0}^{1}\frac{E\big(W_i(u) W_j(1) \big)}{u}du- \int\limits_{0}^{1} \frac{E\big(W_i(1)W_j(v)\big)}{v}dv\big]\\ \nonumber
&&\!\!\!\!\!\!\!\!\!\!\!\!\!\!  =  \nu\gamma_i \gamma_{j-d+1}\big[R_{i,j-d+1}(1,\beta)+\int\limits_{0}^{1} \int\limits_{0}^{1}\frac{R_{i,j-d+1}(u,\beta v)}{uv} dudv-\int\limits_{0}^{1} \frac{R_{i,j-d+1}(u,\beta)}{u}du -\int\limits_{0}^{1} \frac{R_{i,j-d+1}(1,\beta v)}{v}dv\big].%\\=:\nu\gamma_i \gamma_{j-d+1} \big[R_{i,j-d+1}(1,\beta)+I_1-I_2-I_3 \big].
 \end{eqnarray}
 Observe that by two changes of variables and the homogeneity of order 1 of $R_{i,j-d+1}$:
 \begin{eqnarray*}
 &&\quad \int\limits_{0}^{1} \int\limits_{0}^{1}\frac{R_{i,j-d+1}(u,\beta v)}{uv} dudv=\int\limits_{0}^{1} \int\limits_{0}^{v}\frac{R_{i,j-d+1}(u,\beta v)}{uv} dudv+\int\limits_{0}^{1} \int\limits_{0}^{u}\frac{R_{i,j-d+1}(u, \beta v)}{uv} dvdu\\
 &&=\int\limits_{0}^{1} \int\limits_{0}^{1}\frac{R_{i,j-d+1}(vu,\beta v)}{uv} dudv+\int\limits_{0}^{1} \int\limits_{0}^{1}\frac{R_{i,j-d+1}(u, \beta vu)}{uv} dvdu\\
 &&=\int\limits_{0}^{1} \frac{R_{i,j-d+1}(u,\beta)}{u}du +\int\limits_{0}^{1} \frac{R_{i,j-d+1}(1,\beta v)}{v}dv.
  \end{eqnarray*}
  Hence the covariance is equal to $ \nu\gamma_i \gamma_{j-d+1}R_{i,j-d+1}(1,\beta)$.
   \hfill$\Box$

\medskip

\noindent \textbf{Proof of Theorem \ref{theo2}}
  From the uniform consistency of the tail copula estimators and the continuity of the tail copulas we have $\hat{H}_{1j}^{-1} \xrightarrow{P} H_{1j}^{-1},\, j=1, \ldots, d$.  This in combination with  (\ref{ahm}) and Corollary \ref{cor2} yields
   %  Theorem 7.2.1 in \cite{de2006extreme} yields to
%$$ \sup_{0 \leqslant x,y \leqslant t} \lvert \hat{R}(x,y)- R(x,y) \rvert \xrightarrow{p} 0,$$
%which yields
%$\tilde{H}^{-1} \xrightarrow{P} H^{-1}.$
%Then
\begin{equation}
     \sqrt{k}(\hat{\gamma}_{1,d}-\gamma_1)=\sqrt{k}(\hat{\gamma}_1-\gamma_1)+\sum_{j=2}^d \frac{{\gamma}_{1}}{{\gamma}_{j}} \frac{H_{1j}^{-1}}{H_{11}^{-1}}\sqrt{k}(\hat{\gamma}_{j+}-\hat{\gamma}_{j})
     %+\frac{{\gamma}_{1}}{{\gamma}_{3}} \frac{H_{13}^{-1}}{H_{11}^{-1}}\sqrt{k}(\hat{\gamma}_{3+}-\hat{\gamma}{_3})+\ldots+\frac{{\gamma}_{1}}{{\gamma}_{d}} % \frac{H_{1d}^{-1}}{H_{11}^{-1}}\sqrt{k}(\hat{\gamma}_{d+}-\hat{\gamma}{_d})
     +o_p(1). \label{eq8}
     \end{equation}
Now Corollary \ref{cor2} and the  continuous mapping theorem yield  (\ref{eq5_}). \hfill$\Box$

\medskip

\noindent \textbf{Remark 5} The assumption of invertibility of $H$ can be weakened somewhat. In the bivariate case, this weakened assumption follows already from the other assumptions.  Therefore the invertibility of $H$ is not assumed there.

\bibliographystyle{chicago}
\bibliography{References}
\end{document}